\documentclass[12pt,twoside]{amsart}
\usepackage{amssymb}
\usepackage{amscd}

\title[Koll\'ar's injectivity theorem]
{A transcendental approach to Koll\'ar's injectivity theorem}
\author[Osamu Fujino]{Osamu Fujino}
\keywords{multiplier ideal sheaves, vanishing 
theorems, harmonic forms, $\bar\partial$-equations}
\subjclass[2010]{Primary 32L10; Secondary 32W05.}
\date{2011/4/30, version 1.25}
\address{Department of Mathematics, Faculty of Science, 
Kyoto University, Kyoto 606-8502, Japan}  
\email{fujino@math.kyoto-u.ac.jp}
\newcommand{\lla}[0]{{\langle\!\langle}}
\newcommand{\rra}[0]{{\rangle\!\rangle}}
\newcommand{\Hom}[0]{{\operatorname{Hom}}}
\newcommand{\xid}[0]{{\operatorname{id}}}
\newcommand{\Pic}[0]{{\operatorname{Pic}}}
\newcommand{\loc}[0]{{\operatorname{loc}}}
\newcommand{\codim}[0]{{\operatorname{codim}}}
\newcommand{\Ker}[0]{{\operatorname{Ker}}}
\newcommand{\xId}[0]{{\operatorname{Id}}}
\newcommand{\Nak}[0]{{\operatorname{Nak}}}

\newcommand{\xIm}[0]{{\operatorname{Im}}}
\newcommand{\Dom}[0]{{\operatorname{Dom}}}
\newcommand{\xdiv}[0]{{\operatorname{div}}}

\newtheorem{thm}{Theorem}[section]
\newtheorem{lem}[thm]{Lemma}
\newtheorem{cor}[thm]{Corollary}
\newtheorem{prop}[thm]{Proposition}
\newtheorem{claim}{Claim}

\theoremstyle{definition}
\newtheorem{ex}[thm]{Example}
\newtheorem{defn}[thm]{Definition}
\newtheorem{rem}[thm]{Remark}
\newtheorem*{ack}{Acknowledgments}       
\newtheorem{say}[thm]{}
\begin{document}
\bibliographystyle{amsalpha+}
\begin{abstract}
We treat Koll\'ar's injectivity theorem from the analytic 
(or differential geometric) viewpoint. 
More precisely, we give a curvature condition which implies 
Koll\'ar type cohomology injectivity theorems. 
Our main theorem is formulated for a compact K\"ahler manifold, but 
the proof uses the space of harmonic forms on a Zariski open set with a 
suitable complete K\"ahler metric. 
We need neither covering tricks, desingularizations, nor Leray's 
spectral sequence. 
\end{abstract}
\maketitle
\tableofcontents

\section{Introduction}
In \cite{k1}, J\'anos Koll\'ar proved the following theorem. 
We call it {\em{Koll\'ar's original injectivity theorem}} in this paper. 

\begin{thm}[{cf.~\cite[Theorem 2.2]{k1}}]\label{th0} 
Let $X$ be a smooth projective variety defined over 
an algebraically closed field of characteristic 
zero and let $L$ be a semi-ample line bundle on $X$. 
Let $s$ be a nonzero holomorphic section of $L^{\otimes k}$ 
for some $k>0$. Then 
$$
\times s:H^q(X, K_X\otimes L^{\otimes m})\to 
H^q(X, K_X\otimes L^{\otimes m+k}) 
$$ 
is injective for every $q\geq 0$ and every $m\geq 1$, where $K_X$ 
is the canonical line bundle of $X$. 
Note that $\times s$ is the homomorphism induced 
by the tensor product with $s$. 
\end{thm}
The following theorem is the main result of this paper. 
It is an analytic formulation of Koll\'ar type cohomology injectivity theorem. 

\begin{thm}[Main Theorem]\label{main}
Let $X$ be an $n$-dimensional compact K\"ahler manifold. 
Let $(E,h_E)$ $($resp.~$(L,h_L)$$)$ be a holomorphic vector $($resp.~line$)$ 
bundle on $X$ with a smooth hermitian metric $h_E$ $($resp.~$h_L$$)$. 
Let $F$ be a holomorphic line bundle on $X$ with a singular hermitian 
metric $h_F$. 
Assume the following conditions. 
\begin{itemize}
\item[(i)] There exists a subvariety $Z$ of $X$ such that 
$h_F$ is smooth on $X\setminus Z$. 
\item[(ii)] $\sqrt{-1}\Theta (F)\geq -\gamma$ in the sense of 
currents, where $\gamma$ is a smooth 
$(1,1)$-form on $X$. 
\item[(iii)] $\sqrt{-1}(\Theta (E)+\xId_E\otimes \Theta
(F))\geq _{\Nak}0$ on $X\setminus Z$. 
\item[(iv)] $\sqrt {-1} (\Theta (E)+\xId_E \otimes \Theta(F)
-\varepsilon \xId _E\otimes \Theta (L))\geq _{\Nak}0$ on 
$X\setminus Z$ for some positive constant $\varepsilon$. 
\end{itemize}
Here, $\geq _{\Nak}0 $ means the Nakano semi-positivity. 
Let $s$ be a nonzero holomorphic section of $L$. 
Then the multiplication homomorphism 
$$
\times s : 
H^q(X, K_X\otimes E\otimes F\otimes \mathcal J(h_F))\to 
H^q(X, K_X\otimes E\otimes F\otimes \mathcal J(h_F)\otimes L)
$$ 
is injective for every $q\geq 0$, where $\mathcal J(h_F)$ is the 
multiplier ideal sheaf associated to the singular hermitian metric 
$h_F$ of $F$. 
\end{thm}

The formulation of Theorem \ref{main} 
was inspired by Ohsawa's injectivity 
theorem (see \cite{ohsawa}).  
Although the assumptions in Theorem \ref{main} may look artificial 
for algebraic geometers, 
our main theorem 
is useful and have potentiality for various generalizations. 
As a direct consequence of Theorem \ref{main}, 
we have the following corollary. 

\begin{cor}\label{main-cor} 
Let $X$ be an $n$-dimensional compact K\"ahler manifold. 
Let $(E,h_E)$ $($resp.~$(L,h_L)$$)$ be a holomorphic vector $($resp.~line$)$ 
bundle on $X$ with a smooth hermitian metric $h_E$ $($resp.~$h_L$$)$. 
Let $F$ be a holomorphic line bundle on $X$. 
Assume the following conditions. 
\begin{itemize}
\item[(a)] There exists an effective Cartier divisor $D$ on $X$ such that 
$\mathcal O_X(D)\simeq F^{\otimes k}$ for some positive integer $k$. 
\item[(b)] $\sqrt{-1}\Theta (E)
\geq _{\Nak}0$. 
\item[(c)] $\sqrt {-1} (\Theta (E)
-\varepsilon \xId _E\otimes \Theta (L))\geq _{\Nak}0$ 
for some positive constant $\varepsilon$. 
\end{itemize}
Let $s$ be a nonzero holomorphic section of $L$. 
Then the multiplication homomorphism 
$$
\times s : 
H^q(X, K_X\otimes E\otimes F\otimes \mathcal J)\to 
H^q(X, K_X\otimes E\otimes F\otimes \mathcal J\otimes L)
$$ 
is injective for every $q\geq 0$, where $\mathcal J=\mathcal J(\frac{1}{k}D)$ is the 
multiplier ideal sheaf associated 
to $\frac{1}{k}D$ {\em{(}}cf.~{\em{Definition \ref{28}}}{\em{)}}. 
\end{cor}

One of the advantages of our formulation is that we are released from  
sophisticated algebraic geometric methods 
such as desingularizations, 
covering tricks, Leray's spectral sequence, and so on both 
in the 
proof and in 
various applications (see, for example, the proof of Proposition \ref{cor03}). 
The main ingredient of our proof of Theorem \ref{main} 
is Nakano's identity (see Proposition \ref{naka}). 

We note that 
there are many contributors (Koll\'ar, Esnault--Viehweg, 
Kawamata, Ambro, ...) to this kind of cohomology injectivity 
theorem. 
We just mention that the first result was obtained 
by Tankeev \cite[Proposition 1]{tan}. 
It inspired Koll\'ar to obtain his famous 
injectivity theorem (see \cite{k1} or Theorem \ref{th0}). 
After \cite{k1}, many generalizations of Theorem \ref{th0} 
were obtained (see the books \cite{ev} and \cite{k2}). 
Koll\'ar did not refer to \cite{enoki} in \cite{k2}. 
However, we think that \cite{enoki} is the first paper 
where Koll\'ar's injectivity theorem is proved (and generalized) 
by differential 
geometric arguments. 

Let us recall Enoki's theorem \cite[Theorem 0.2]{enoki}, 
which is a very special case of Theorem \ref{main}, 
for the reader's convenience. 
To recover Corollary \ref{eno} from Theorem \ref{main}, it is sufficient to 
put $E=\mathcal O_X$, $F=L^{\otimes m}$, and $L=L^{\otimes k}$. 
The reader who reads Japanese can find \cite{fujino3} 
useful. It is a survey on Enoki's injectivity theorem. 

\begin{cor}[Enoki]\label{eno}
Let $X$ be an $n$-dimensional 
compact K\"ahler manifold and 
let $L$ be a semi-positive holomorphic line bundle on $X$. 
Suppose $L^{\otimes k}$, $k>0$, admits a nonzero 
global holomorphic section $s$. 
Then 
$$
\times s: H^q(X, K_X\otimes L^{\otimes m})\to 
H^q(X, K_X\otimes L^{\otimes m+k})
$$ is injective for every $m>0$ and every $q\geq 0$. 
\end{cor}
We recall Enoki's idea of the proof in \cite{enoki} 
because we will use the same idea to 
prove Theorem \ref{main}. 

\begin{say}[Enoki's proof] 
From now on, we assume that $k=m=1$ for simplicity. 
It is well known that the cohomology group 
$H^q(X, K_X\otimes L^{\otimes l})$ 
is represented by 
the space of harmonic forms 
$\mathcal H^{n,q}(L^{\otimes l})=\{u :{\text{smooth 
$L^{\otimes l}$-valued 
$(n,q)$-form on $X$ such that }} \bar\partial u=0, 
D''^*_{L^{\otimes l}}u=0\}$, where $D''^*_{L^{\otimes l}}$ is 
the formal adjoint of $\bar\partial$. 
We take 
$u\in \mathcal H^{n,q}(L)$. Then, $\bar\partial (su)=0$ because 
$s$ is holomorphic. We can check that $D''^*_{L^{\otimes 2}}(su)=0$ by 
using Nakano's identity and the semi-positivity of $L$. 
Thus, $s$ induces $\times s: \mathcal H^{n,q}(L)\to  
\mathcal H^{n,q}(L^{\otimes 2})$. 
Therefore, the required injectivity is 
obvious. 
\end{say}

Enoki's theorem contains Koll\'ar's original injectivity theorem 
(cf.~Theorem \ref{th0}) 
by the following 
well-known lemma. 

\begin{lem}
Let $L$ be a semi-ample line bundle on a smooth projective manifold $X$. 
Then $L$ is semi-positive. 
\end{lem}
\begin{proof}
There exists a morphism $f=\Phi_{|L^{\otimes m}|}:X\to 
\mathbb P^N$ induced by the 
complete linear system $|L^{\otimes m}|$ for some 
$m>0$ because $L$ is semi-ample. 
Let $h$ 
be a smooth hermitian metric on 
$\mathcal O_{\mathbb P^N}(1)$ with positive definite curvature. 
Then $(f^*h)^{\frac{1}{m}}$ 
is a smooth hermitian metric on $L$ whose 
curvature is semi-positive. 
\end{proof}

\begin{rem}
Let $X$ be a complex analytic space and 
let $\mathcal E$ be a coherent sheaf on $X$. 
In order to prove $H^p(X, \mathcal E)=0$, it is 
sufficient to construct a homomorphism 
$\varphi:\mathcal E\to \mathcal F$ of coherent sheaves on $X$ such 
that 
the induced map $H^p(X,\mathcal E)\to H^p(X, \mathcal F)$ is 
{\em{injective}} and that $H^p(X, \mathcal F)=0$. 
This simple observation plays crucial roles 
for various vanishing theorems on toric varieties (see, 
for example, \cite{fujino-multi} and \cite{fujino-vani}). 
Anyway, injectivity theorems sometimes 
are very useful in proving various vanishing theorems. 
See the proof of Corollary \ref{cor2} below. 
\end{rem}

We quickly review Koll\'ar's proof of his injectivity theorem in \cite{k2}, 
which is much simpler than Koll\'ar's original proof in \cite{k1}, for the 
reader's convenience. 

\begin{say}[Koll\'ar's proof]\label{15} 
Let $X$ be a smooth projective $n$-fold 
and let $L$ be a (not necessarily semi-ample) line bundle on $X$. 
Let $s$ be a non-zero holomorphic 
section of $L^{\otimes 2}$. 
Assume that $D=(s=0)$ is a {\em{smooth}} divisor on $X$ for simplicity. 
We can take a double cover $\pi:Z\to X$ ramifying along 
$D$. 
By the Hodge decomposition, we obtain a surjection $H^q(Z, \mathbb C_Z)
\to H^q(X, \mathcal O_Z)$ for every $q$. 
By taking the anti-invariant part of the covering involution, we 
obtain that $H^q(X, G)\to H^q(X, L^{-1})$ is surjective for every $q$, 
where $\pi_*\mathbb C_Z=\mathbb C_X\oplus G$ is the eigen-sheaf decomposition. 
It is not difficult to see that there exists a factorization 
$H^q(X, G)\to H^q(X, L^{-1}\otimes \mathcal O_X(-D))
\to H^q(X, L^{-1})$ for every $q$. 
Therefore, $\times s:H^q(X, K_X\otimes L)\to 
H^q(X, K_X\otimes L\otimes \mathcal O_X(D))$ is injective by the Serre 
duality. 
In general, $D$ is not necessarily {\em{smooth}}. 
So, we have to use sophisticated algebraic geometric methods such as 
desingularizations, relative vanishing theorems, Leray's spectral sequences, 
and so on, even when $X$ is smooth and $L$ is free. 
\end{say} 

\begin{rem}
As we saw in \ref{15}, thanks to the Serre duality, the injectivity of 
$H^q(X, K_X\otimes L)\to H^q(X, K_X\otimes L\otimes \mathcal 
O_X(D))$ is equivalent to the surjectivity of $H^{n-q}(X, L^{-1}\otimes 
\mathcal O_X(-D))\to H^{n-q}(X, L^{-1})$. 
However, injectivity seems to be much better and more natural for 
some applications and generalizations. See Section \ref{sec4}. 
\end{rem}

Roughly speaking, Koll\'ar's geometric 
proof in \cite{k2} (and Esnault--Viehweg's proof 
in \cite{ev}) depends on the Hodge decomposition, or the degeneration 
of the Hodge to de Rham type 
spectral sequence. So, it works only when 
$E$ is a unitary flat vector bundle (see \cite[9.17 Remark]{k2}). 
On the other hand, our analytic proof (and 
the proofs in \cite{enoki}, \cite{ohsawa}, and \cite{take}) 
relies on the harmonic representation of the 
cohomology groups. We do not know the true relationship between 
the geometric proof and the analytic one. 

\begin{say}[More advanced topics] 
In \cite{fujino2}, we prove a relative version of Theorem \ref{main}. 
In that case, $X$ is not necessarily compact. 
When $X$ is not compact, a {\em{locally}} 
square integrable differential form $u$ on $X$ is 
not necessarily {\em{globally}} square integrable. 
So, we use Ohsawa--Takegoshi's twisted version of 
Nakano's identity to control 
the asymptotic behavior of the $L^2$-norm of $u$ around 
the boundary of $X$. Thus, we need much more 
analytic methods for the relative setting. 

In \cite[Chapter 2]{fujino-book}, 
\cite{fujino-on}, and \cite[Sections 5 and 6]{fujino-fundamental}, 
we develop the geometric approach (see \ref{15}) to 
obtain a very important generalization 
of Koll\'ar's injectivity theorem. In those papers, 
we consider mixed Hodge structures on compact support 
cohomology groups. 
Roughly speaking, 
the decomposition 
$$
H^n_c(X\setminus \Sigma, \mathbb C)\simeq \bigoplus _{p+q=n}
H^q(X, \Omega^p_X(\log \Sigma)\otimes \mathcal O_X(-\Sigma))
$$ 
where $X$ is a smooth projective variety and 
$\Sigma$ is a simple normal crossing divisor on $X$ 
produces a generalization of Koll\'ar type cohomology injectivity theorem. 
The reader can find a thorough treatment of our geometric 
approach in \cite[Chapter 2]{fujino-book}. 
We have already obtained many applications 
for the log minimal model program in \cite{fujino-qlog}, \cite{fujino-book}, 
\cite{fujino-eff1}, \cite{fujino-eff2}, 
\cite{fujino-theory}, \cite{fst}, 
\cite{fujino-non}, \cite{fujino-fundamental}, and \cite{fujino-surface}. 

By our experience, we know that 
Koll\'ar type injectivity theorems play crucial 
roles for the study of base point free theorems and 
the abundance conjecture for log canonical 
pairs 
(cf.~\cite{fukuda}, \cite{fujino-bpf}, and so on). 
\end{say}

We summarize the contents of this paper. 
In Section \ref{sec2}, we fix notation and collect basic results. 
Section \ref{sec3} is the proof of the main theorem:~Theorem \ref{main}. 
We will represent the cohomology groups by the spaces of harmonic forms 
on a Zariski open set with a suitable complete K\"ahler metric. 
We will use $L^2$-estimates for $\bar\partial$-equations on 
complete K\"ahler manifolds (see Lemma \ref{3.2}). It is 
a key point of our proof.  
In Section \ref{sec4}, we treat Koll\'ar type injectivity 
theorem, Esnault--Viehweg type injectivity theorem, 
and Kawamata--Viehweg--Nadel type vanishing 
theorem as applications of Theorem \ref{main}. 
We recommend the reader to compare 
them with usual algebraic geometric ones. 
We note that we discuss them in a more general relative 
setting in \cite{fujino2}. 

\begin{ack} 
The first version of this paper was written 
in Nagoya in 2006 and was circulated as arXiv:0704.0073. 
The author would like to thank Professor Takeo Ohsawa for answering 
his questions. 
He was partially supported by The Sumitomo Foundation 
and by the Grant-in-Aid for Young Scientists (A) 
$\sharp$17684001 from 
JSPS when he prepared the first version. 
He thanks Doctor Dano Kim for useful comments. 
He revised this paper in Kyoto in 2010. 
He was partially supported by The Inamori Foundation and by the 
Grant-in-Aid for Young Scientists (A) $\sharp$20684001 from JSPS. 
He thanks the referee for useful comments and informing him of 
the papers \cite{siu} and \cite{ep}. 
\end{ack}

\section{Preliminaries}\label{sec2}
In this section, we collect basic definitions and results in algebraic and 
analytic geometries. For details, see, for example, \cite{dem}. 

\begin{say}[Singular hermitian metric] 

Let $L$ be a holomorphic line bundle on a complex manifold $X$. 

\begin{defn}[Singular hermitian metric] 
A {\em{singular hermitian metric}} on $L$ is a metric which is given 
in every trivialization $\theta:L|_{\Omega}\simeq \Omega \times 
\mathbb C$ by 
$$
\| \xi\|=|\theta (\xi)|e^{-\varphi(x)},\ \ \  x\in \Omega, \ \xi \in L_x, 
$$
where $\varphi\in L^1_{\loc}(\Omega)$ is an arbitrary function, 
called the {\em{weight}} of the metric with respect to the trivialization $\theta$. 
Here, $L^1_{\loc}(\Omega)$ is the space of the locally integrable functions 
on $\Omega$. 
\end{defn}

The following singular hermitian metrics play important roles 
in the study of higher dimensional algebraic varieties. 

\begin{ex}\label{ex1}
Let $D=\sum \alpha _j D_j$ be a divisor with coefficients $\alpha _j\in \mathbb N$. 
Then $\mathcal O_X(D)$ is equipped with a natural singular hermitian metric 
as follows. 
Let $f$ be a local section of $\mathcal O_X(D)$, viewed as a meromorphic 
function such that $\xdiv (f)+D\geq 0$. 
We define 
$\|f\|^2=|f|^2\in [0,\infty]$. If $g_j$ is a generator of the ideal 
of $D_j$ on an open set $\Omega \subset X$, then 
the weight corresponding to this metric is $\varphi=\sum _j \alpha _j \log 
|g_j|$. 
It is obvious that this metric is a smooth hermitian metric on $X\setminus D$ 
and its curvature is zero on $X\setminus D$. 
Let $L$ be a holomorphic line bundle on $X$. 
Assume that $L^{\otimes k}\simeq M\otimes \mathcal O_X(D)$ for 
some holomorphic line bundle $M$ and an effective divisor $D$ on $X$. 
As above, 
$\mathcal O_X(D)$ is equipped with a natural 
singular hermitian metric $h_D$. 
Let $h_M$ be any smooth hermitian metric on $M$. 
Then $L$ has a singular hermitian metric $h_L:=h^{\frac{1}{k}}_Mh^{\frac{1}{k}}_D$. 
Note that $h_L$ is smooth outside $D$ and $\Theta _{h_L}(L)
=\frac{1}{k}\Theta _{h_M}(M)$ on $X\setminus D$. 
\end{ex}
\end{say}

\begin{say}[Multiplier ideal sheaf] 
The notion of multiplier ideal sheaves introduced by Nadel \cite{nadel} 
is very important in recent developments of 
complex and algebraic geometries (cf.~\cite[Part Three]{laz}). 

\begin{defn}[(Quasi-)plurisubharmonic function and 
multiplier ideal sheaf] 
A function $u:\Omega\to [-\infty, \infty)$ defined on an open set 
$\Omega\subset \mathbb C^n$ is called {\em{plurisubharmonic}} 
(psh, for short) if 
\begin{itemize}
\item[1.] $u$ is upper semi-continuous, and 
\item[2.] for every complex line $L\subset \mathbb C^n$, 
$u|_{\Omega\cap L}$ is subharmonic on $\Omega\cap L$, that is, 
for every $a\in \Omega$ and $\xi \in \mathbb C^n$ 
satisfying $|\xi|<d(a, \Omega^c)$, the function $u$ satisfies 
the mean inequality 
$$
u(a)\leq \frac{1}{2\pi} \int ^{2\pi}_{0} u(a+e^{i\theta}\xi)d\theta. 
$$
\end{itemize} 
Let $X$ be an $n$-dimensional complex manifold. 
A function $\varphi:X\to [-\infty, \infty)$ is said to be 
{\em{plurisubharmonic}} (psh, for short) if 
there exists an open cover $X=\bigcup _{i\in I}U_i$ such that 
$\varphi|_{U_i}$ is plurisubharmonic on $U_i$ 
($\subset \mathbb C^n$) for every $i$. 
A {\em{smooth strictly plurisubharmonic function}} $\psi$ on $X$ is a 
smooth function on $X$ such that $\sqrt{-1} \partial \bar\partial \psi$ is 
a positive definite smooth $(1,1)$-form. 
A {\em{quasi-plurisubharmonic}} 
(quasi-psh, for short) function is a function $\varphi$ which 
is locally equal to the sum of a psh function and of a smooth function. 
If $\varphi$ is a quasi-psh function on a complex manifold $X$, 
the {\em{multiplier ideal sheaf}} 
$\mathcal J(\varphi)\subset \mathcal O_X$ is 
defined by 
$$
\Gamma (U, \mathcal J(\varphi))=\{f\in \mathcal O_X(U);\ |f|^2e^{-2\varphi}\in 
L^1_{\loc}(U)\} 
$$ 
for every open set $U\subset X$. 
Then it is known that 
$\mathcal J(\varphi)$ is a coherent ideal sheaf of $\mathcal O_X$. 
See, for example, \cite[(5.7) Proposition]{dem}. 
\end{defn}

\begin{rem}
By the assumption (ii) in Theorem \ref{main}, 
the weight of the singular hermitian metric $h_F$ is a quasi-psh function 
on every trivialization. 
So, we can define multiplier ideal sheaves locally and 
check that they are independent of trivializations. 
Thus, we can define 
the multiplier ideal sheaf globally and denote it by $\mathcal J(h_F)$, 
which is an abuse of notation.  
It is a coherent ideal 
sheaf on $X$. 
\end{rem}

\begin{ex}
Let $X=\{\, z\in \mathbb C \, | \, |z|<r \}$ for some $0<r<1$ and 
let $L$ be a trivial line bundle on $X$. 
We consider a singular hermitian metric $h_L=\exp(\sqrt{-\log |z|^2})$ of $L$. 
Then $h_L$ is smooth outside the origin $0\in X$. The weight of $h_L$ is $\varphi
=-\frac{1}{2}\sqrt{-\log |z|^2}$ and $\varphi$ is a psh function on $X$. 
The Lelong number of $\varphi$ at $0$ is 
$$\liminf_{z\to 0}\frac{\varphi(z)}{\log |z|}
=0. 
$$ 
Thus, we have $\mathcal J(h_L)\simeq \mathcal O_X$ by Skoda. 
Note that $\varphi$ is smooth outside $0$, which is an analytic 
subvariety of $X$. However, $\varphi$ does not have analytic singularities 
around $0$.  
\end{ex}

\begin{defn}\label{28}
Let $X$ be a complex manifold and let $D=\sum \alpha_j D_j$ be an effective 
$\mathbb Q$-divisor on $X$. 
Let $g_j$ be a generator of the ideal of $D_j$ on an open set 
$\Omega\subset X$. 
We put $\mathcal J(D):=\mathcal J(\varphi)$, 
where $\varphi=\sum _j\alpha_j\log |g_j|$. 
Since $\mathcal J(\varphi)$ is independent of the choice of the generators 
$g_j$'s, $\mathcal J(D)$ is a well-defined coherent ideal 
sheaf on $X$. We call $\mathcal J(D)$ 
the multiplier ideal sheaf associated to the effective $\mathbb Q$-divisor $D$. 
We say that the divisor $D$ is {\em{integrable}} 
at a point $x_0\in X$ if the function 
$\prod |g_j|^{-2\alpha_j}$ 
is integrable on a neighborhood of 
$x_0$, equivalently, $\mathcal J(D)_{x_0}=\mathcal O_{X, x_0}$. 
Let $D'$ be another effective 
$\mathbb Q$-divisor on $X$. 
Then, 
$\mathcal J(D)=\mathcal J(D+\varepsilon D')$ for 
$0<\varepsilon \ll 1$, $\varepsilon \in \mathbb Q$. 
\end{defn}

\begin{rem}
In Definition \ref{28}, 
$D$ is integrable at $x_0$ if and only if 
the pair $(X, D)$ is {\em{Kawamata log terminal}} 
({\em{klt}}, for short) in a neighborhood of $x_0$ (cf.~\cite[Definition 2.34]{km}). 
\end{rem}

\begin{ex}
Let $h_L$ be the singular hermitian metric defined in Example \ref{ex1}. 
Then the weight of the singular hermitian metric $h_L$ is a quasi-psh function 
on every trivialization. 
Therefore, the multiplier ideal sheaf $\mathcal J(h_L)$ is well-defined 
and $\mathcal J(h_L)=\mathcal J(\frac{1}{k}D)$.  
\end{ex}
\end{say}
\begin{say}[Hermitian and K\"ahler geometries]
We collect the basic notion and results 
of hermitian and K\"ahler geometries (see also \cite{dem}). 

\begin{defn}[Chern connection 
and its curvature form]\label{conecone}
Let $X$ be a complex hermitian manifold and 
let $(E,h)$ be a holomorphic hermitian vector bundle on $X$. 
Then there exists the {\em{Chern connection}} $D=D_{(E, h)}$, 
which can be split in a unique way as a sum of a $(1,0)$ and of 
a $(0,1)$-connection, $D=D'_{(E, h)}+D''_{(E,h)}$. 
By the definition of the Chern connection, 
$D''=D''_{(E,h)}=\bar\partial$. 
We obtain the {\em{curvature form}} $\Theta(E)=\Theta_{(E, h)}
=\Theta_h:=D^2_{(E, h)}$. The subscripts might be suppressed 
if there is no danger of confusion. 
\end{defn}

Let $U$ be a small open set of $X$ and let $(e_\lambda)$ be a local 
holomorphic frame of $E|_{U}$. 
Then the hermitian metric $h$ is given by the hermitian 
matrix $H=(h_{\lambda\mu})$, 
$h_{\lambda\mu}=h(e_\lambda, e_\mu)$, on $U$. 
We have $h(u, v)={}^tuH\bar v$ on $U$ for 
smooth sections $u$, $v$ of $E|_{U}$. 
This implies that $h(u,v)=\sum_{\lambda, \mu}u_\lambda h_{\lambda\mu}
\bar{v}_{\mu}$ for $u=\sum e_iu_i$ and $v=\sum e_j v_j$. 
Then we obtain that 
$\sqrt{-1}\Theta_h(E)=\sqrt{-1}\bar\partial(\overline{H}^{-1}\partial \overline{H})$ 
and ${}^t\overline{(\sqrt{-1}{}^t\Theta_h(E)H)}=\sqrt{-1}{}^t\Theta_h(E)H$ on 
$U$. 

\begin{defn}[Inner product]\label{d216}
Let $X$ be an $n$-dimensional complex 
manifold with the hermitian metric $g$. 
We denote by $\omega$ the {\em{fundamental form}} of $g$. 
Let $(E,h)$ be a hermitian vector bundle on $X$, 
and $u,v$ 
are $E$-valued $(p,q)$-forms with measurable coefficients, 
we set 
$$
\|u\|^2=\int _X |u|^2dV_\omega, 
\ \lla u, v \rra= \int _X \langle u, v \rangle 
dV_{\omega}, 
$$ 
where $|u|$ is the pointwise norm induced by $g$ and 
$h$ on $\Lambda^{p,q}T^*_X\otimes E$, and 
$dV_\omega=\frac{1}{n!}\omega^n$. 
More explicitly, $\langle u, v\rangle dV_\omega
={}^t u\wedge H\overline{\ast v}$, where $\ast$ 
is the {\em{Hodge star operator}} relative 
to $\omega$ and $H$ is the (local) matrix representation of $h$. 
When we need to emphasize the metrics, we write 
$|u|_{g,h}$, and so on. 
\end{defn}

Let $L^{p,q}_{(2)}(X,E)(=
L^{p.q}_{(2)}(X,(E,h)))$ 
be the space of square integrable $E$-valued $(p,q)$-forms on $X$. 
The inner product was defined in 
Definition \ref{d216}. 
When we emphasize the metrics, 
we write $L^{p,q}_{(2)}(X,E)_{g,h}$, where 
$g$ (resp.~$h$) is the hermitian metric of $X$ 
(resp.~$E$). 
As usual one can view $D'$ and $D''$ as closed and densely 
defined operators on the Hilbert space $L^{p,q}_{(2)}(X, E)$. 
The formal adjoints 
${D'}^*$, ${D''}^*$ also have closed extensions in 
the sense of distributions, which do not necessarily 
coincide with the Hilbert space adjoints in the sense of 
Von Neumann, since 
the latter ones may have strictly smaller domains. 
It is well known, however, that the domains coincide if 
the hermitian metric of $X$ is complete. See Lemma \ref{density} below. 

\begin{defn}[Nakano positivity and semi-positivity] 
Let $(E, h)$ be a holomorphic vector bundle on a complex manifold $X$ with 
a smooth hermitian metric $h$. 
Let $\Xi$ be a $\Hom (E, E)$-valued $(1,1)$-form such that 
${}^t\overline{({}^t\Xi h)}={}^t\Xi h$. 
Then $\Xi$ is said to be {\em{Nakano positive}} (resp.~{\em{Nakano 
semi-positive}}) 
if the hermitian form on $T_X\otimes E$ associated to 
${}^t\Xi h$ is positive definite 
(resp.~semi-definite). 
We write $\Xi>_{\Nak}0$ (resp.~$\geq_{\Nak}0$). We note that 
$\Xi_1>_{\Nak}\Xi_2$ (resp.~$\Xi_1\geq _{\Nak}\Xi_2$) means that 
$\Xi_1-\Xi_2>_{\Nak} 0$ (resp.~$\geq _{\Nak} 0$). 
A holomorphic vector bundle $(E, h)$ is said to be {\em{Nakano positive}} 
(resp.~{\em{Nakano 
semi-positive}}) if $\sqrt{-1}\Theta(E)>_{\Nak}0$ (resp.~$\geq _{\Nak}0$). 
We usually omit \lq\lq Nakano\rq\rq when $E$ is a line bundle. 
\end{defn}

\begin{defn}[Graded Lie bracket] 
Let $C^{\infty}(X, \Lambda ^{p,q}T^*_X\otimes E)$ be the 
space of the smooth $E$-valued $(p,q)$-forms on $X$. 
If $A, B$ are the endomorphisms of pure degree of the graded module 
$M^{\bullet}
=C^{\infty}(X, \Lambda ^{\bullet, \bullet}T^*_X\otimes E)$, their 
{\em{graded Lie bracket}} is defined by 
$$[A, B]=AB-(-1)^{\deg A\deg B}BA.$$ 
\end{defn}

Let us recall Nakano's identity, which 
is one of the main ingredients of the 
proof of our main theorem:~Theorem \ref{main}. 

\begin{prop}[Nakano's identity]\label{naka} 
We further assume that $g$ is K\"ahler. Let 
$$\Delta'=D'D'^*+D'^*D'$$ and 
$$\Delta''=D''D''^*+D''^*D''$$ be the complex Laplace operators acting 
on $E$-valued forms. 
Then 
$$
\Delta''=\Delta'+[\sqrt{-1}\Theta (E), \Lambda],  
$$ 
where $\Lambda$ is the adjoint of $\omega\wedge \cdot \ $. 
\end{prop}

The following lemma is now classical. 
See, for example, \cite[Lemme 4.3]{dem-est}.  

\begin{lem}[Density lemma]\label{density}
If $g$ is complete, then $C^{p,q}_0(X,E)$ is dense in $\Dom D''^*\cap \Dom \bar 
\partial$ with respect to the graph norm  
$$
u\mapsto \| u\|+\|\bar\partial u\|+\| D''^*u\|, 
$$ 
where $C^{p,q}_0(X,E)$ is the space of the $E$-valued smooth 
$(p,q)$-forms on $X$ with compact supports and 
$\Dom D''^*$ $($resp.~$\Dom \bar\partial$$)$ is the domain of $D''^*$ 
$($resp.~$\bar\partial$$)$. 
\end{lem}

Combining Proposition \ref{naka} with Lemma \ref{density}, 
we obtain the following formula. 

\begin{prop}\label{tsuika}
Let $u$ be a square integrable $E$-valued $(n,q)$-form 
on $X$ with $\dim X=n$ 
and let $g$ be a complete K\"ahler metric on $X$. 
Let $\omega$ be the fundamental form of $g$. 
Assume that $\sqrt{-1}\Theta(E)
\geq_{\Nak}-c\xId _E\otimes \omega$ for some constant 
$c$. 
Then we obtain that 
$$
\|D''^*u\|^2+
\|\bar\partial u\|^2=\|{D'}^*u \|^2+
\lla\sqrt{-1}
\Theta(E)\Lambda u, u\rra
$$
for every $u\in \Dom D''^*\cap \Dom \bar\partial $. 
\end{prop}

The final remark in this section will play crucial roles in the proof of the main 
theorem:~Theorem \ref{main}. The proof is an easy calculation (cf.~\cite[Lemme 3.3]{dem-est}). 

\begin{rem}\label{rem215}
Let $g'$ be another hermitian metric on $X$ such that 
$g'\geq g$ and $\omega'$ be the fundamental 
form of $g'$. 
Let $u$ be an $E$-valued $(n,q)$-form with measurable coefficients. 
Then, we have 
$|u|^2_{g', h}dV_{\omega'}\leq |u|^2_{g, h}dV_{\omega}$, where 
$|u|_{g',h}$ (resp.~$|u|_{g,h}$) is the pointwise norm induced 
by $g'$ and $h$ (resp.~$g$ and $h$). 
If $u$ is an $E$-valued $(n,0)$-form, then 
$|u|^2_{g', h}dV_{\omega'}= |u|^2_{g, h}dV_{\omega}$. 
In particular, $\|u\|^2$ is independent of $g$ when 
$u$ is an $(n,0)$-form. 
\end{rem}
\end{say}

\section{Proof of the main theorem}\label{sec3}
In this section, we prove the main theorem:~Theorem \ref{main}. 
The idea is very simple. We represent the cohomology groups by the 
space of harmonic  forms on $X\setminus Z$ (not on $X$!). 
The manifold $X\setminus Z$ is not compact. 
However, it is a complete K\"ahler 
manifold and all hermitian metrics are smooth on $X\setminus Z$. 
So, there are no difficulties on $X\setminus Z$. 
Note that we do not 
need the difficult regularization technique for 
quasi-psh functions on K\"ahler manifolds 
(cf.~\cite[Th\'eor\`{e}me 9.1]{dem-est}). 

\begin{proof}[Proof of {\em{Theorem \ref{main}}}] 
Since $X$ is compact, there exists a complete K\"ahler metric $g'$ on $Y:=X
\setminus Z$ such that $g'>g$ on $Y$. 
We sketch the construction of $g'$ because we need some 
special properties of $g'$ in the 
following proof. 
The next lemma is well known. See, for example, \cite[Lemma 5]{dem-coh}. 
\begin{lem}
There exists a quasi-psh function $\psi$ on $X$ such that 
$\psi=-\infty$ on $Z$ with logarithmic poles along $Z$ and $\psi$ 
is smooth outside $Z$. 
\end{lem}
Without loss of generality, we can assume that $\psi< -e$ on $X$. 
We put $\varphi=\frac{1}{\log (-\psi)}$. Then $\varphi$ is a quasi-psh function 
on $X$ and 
$\varphi< 1$. 
Thus, we can take a positive constant $\alpha$ such that 
$\sqrt{-1} \partial \bar\partial 
\varphi+\alpha \omega> 0$ 
on $Y$. Let $g'$ be the K\"ahler metric on $Y$ whose 
fundamental form is $\omega'=\omega+(\sqrt{-1} \partial \bar\partial 
\varphi+\alpha \omega)$. We will show that 
$$
\omega'\geq \partial (\log (\log (-\psi)))\wedge 
\bar\partial (\log (\log (-\psi)))  
$$ if we choose $\alpha\gg 0$. 
We have $$
\bar\partial \varphi =
-\frac{\frac{-\bar\partial \psi}{-\psi}}
{(\log (-\psi))^2}, 
$$ 
and 
\begin{align*}
\partial \bar\partial \varphi &=
2\frac{\frac{-\partial \psi}{-\psi}
\wedge \frac{-\bar\partial \psi}{-\psi}}
{(\log (-\psi))^3}
-\frac{\partial 
(\frac{-\bar\partial \psi}{-\psi})}
{(\log 
(-\psi))^2}\\
&=
2\frac{
\frac{-\partial \psi}{-\psi}
\wedge \frac{-\bar\partial \psi}{-\psi}
}{(\log (-\psi))^3}
-\frac{
\frac{-\partial \bar\partial \psi}{-\psi}
}{(\log (-\psi))^2}
+\frac{
\frac{-\partial \psi \wedge (-\bar\partial \psi)}
{(-\psi)^2}
}{(\log (-\psi))^2}\\ 
&=2\frac{
\frac{-\partial \psi}{-\psi}
\wedge \frac{-\bar\partial \psi}{-\psi}
}{(\log (-\psi))^3}
+
\frac{
\frac{\partial \bar\partial \psi}{-\psi}
}{(\log (-\psi))^2}
+
\frac{
\frac{-\partial \psi \wedge (-\bar\partial \psi)}
{(-\psi)^2}}
{(\log (-\psi))^2}. 
\end{align*}
On the other hand, 
$$
\partial (\log(\log (-\psi)))=\frac{
\frac{-\partial \psi}{-\psi}}{\log (-\psi)}. 
$$ 
Therefore, 
$$
\partial (\log(\log (-\psi)))\wedge \bar \partial (\log(\log (-\psi)))
=\frac{\frac{-\partial \psi \wedge (-\bar\partial \psi)}
{(-\psi)^2}}{(\log (-\psi))^2}. 
$$
This implies 
$$
\omega'\geq \partial (\log (\log (-\psi)))\wedge 
\bar\partial (\log (\log (-\psi)))   
$$ if 
$\alpha \gg 0$. Therefore, $g'$ is a {\em{complete}} K\"ahler metric on $Y$ by 
Hopf--Rinow because 
$\log (\log (-\psi))$ tends to $+\infty$ on $Z$. 
More precisely, $$\eta:=\frac{1}{\sqrt{2}}\log (\log (-\psi))$$ is a smooth 
exhaustive function on $Y$ such that $|d\eta|_{g'}\leq 1$. 
We fix these K\"ahler metrics throughout this proof. 
In general, 
\begin{align*}
L^{n,q}_{(2)}(Y,E\otimes F)&=L^{n,q}_{(2)}(Y,E\otimes F)_{g', h_Eh_F}
\\ &=\overline {\xIm \bar\partial}\oplus 
\mathcal H^{n,q}(E\otimes F)\oplus \overline 
{\xIm D''^*_{E\otimes F}},
\end{align*} 
where 
$$\mathcal H^{n,q}(E\otimes F):=\{u\in L^{n,q}_{(2)}(Y,E\otimes F)\, | \, 
\bar\partial u=D''^*_{E\otimes F}u=0\}$$ is the space of 
the $E\otimes F$-valued harmonic $(n,q)$-forms. 
We note that $u\in \mathcal H^{n,q}(E\otimes F)$ is smooth 
by the regularization theorem for the 
elliptic operator $\Delta''_{E\otimes F}=D''^*_{E\otimes F}\bar\partial 
+\bar\partial D''^*_{E\otimes F}$. 
The claim below is more or less known to 
experts (cf.~\cite[Section 2]{siu}, \cite[Proposition 4.6]{take} and 
\cite[Theorem 4.13]{oh-hon}). 
We write it for the reader's convenience. 

\begin{claim}\label{cla1} 
We have the following equalities and an isomorphism 
of cohomology groups for every $q\geq 0$. 
\begin{align*}
&\overline {\xIm \bar\partial}=\xIm \bar\partial, \ \ \ \overline {\xIm D''^*_{E\otimes 
F}}=\xIm D''^*_{E\otimes F},\ \  \text{and}\\ 
&{{H}}^q(X, K_X\otimes E\otimes F\otimes \mathcal J(h_F))
\simeq \frac{L^{n,q}_{(2)}
(Y, E\otimes F)\cap \Ker \bar \partial}{\xIm \bar\partial}. 
\end{align*}
\end{claim}\label{claim1}
If the claim is true, then ${{H}}^q(X, K_X\otimes E\otimes F\otimes \mathcal J(h_F))
\simeq \mathcal H^{n,q}(E\otimes F)$ because 
$L^{n,q}_{(2)}(Y, E\otimes F)\cap \Ker \bar \partial=\xIm \bar\partial 
\oplus \mathcal H^{n,q}(E\otimes F)$. 

\begin{proof}[Proof of Claim]
First, let $X=\bigcup _{i\in I}U_i$ be a finite Stein cover of $X$ 
such that each $U_i$ is small. 
We can assume that there is 
a small Stein open set $V_i$ of $X$ such that 
$U_i\Subset V_i$ for every $i$ (see the proof of Lemma \ref{3.2}). 
We denote this cover by $\mathcal U=\{U_i\}_{i\in I}$. 
By Cartan and Leray, we obtain 
$$H^q(X, K_X\otimes E\otimes F\otimes \mathcal J(h_F))\simeq 
{\check{H}}^q(\mathcal U, K_X\otimes E\otimes F\otimes \mathcal J(h_F)), $$ 
where the right hand side is the $\check{\rm C}$ech cohomology group calculated 
by $\mathcal U$. 
Let $\{\rho_i\}_{i \in I}$ be a partition of unity 
associated to $\mathcal U$. 
We put $U_{i_0i_1\cdots i_q}=U_{i_0}\cap \cdots \cap U_{i_q}$. 
Then $U_{i_0i_1\cdots i_q}$ is Stein. 
Let $u=\{u_{i_0i_1\cdots i_q}\}$ such that 
$u_{i_0i_1\cdots i_q}\in \Gamma(U_{i_0i_1\cdots i_q}, K_X\otimes E\otimes 
F\otimes \mathcal J(h_F))$ and $\delta u=0$, 
where $\delta$ is the coboundary operator of $\check{\rm C}$ech complexes. 
We put $u^1=\{u^1_{i_0\cdots i_{q-1}}\}$ with 
$u^1_{i_0\cdots i_{q-1}}=\sum _i \rho_i u_{ii_0\cdots i_{q-1}}$. 
Then $\delta u^1=u$ and $\delta (\bar\partial u^1)=0$. 
Thus, we can construct $u^2$ such that $\delta u^2=\bar\partial 
u^1$ as above by using $\{\rho_i\}$. 
By repeating this process, we obtain $\bar\partial u^q\in 
L^{n,q}_{(2)}(Y, E\otimes F)\cap \Ker \bar\partial$ by Remark 
\ref{rem215} because 
$X$ is compact. 
By the standard diagram chasing, we have a homomorphism 
$$
\bar{\alpha}: {\check{H}}^q({\mathcal U}, K_X\otimes 
E\otimes F\otimes \mathcal J(h_F))
\to 
\frac{L^{n,q}_{(2)}(Y, E\otimes F)
\cap \Ker \bar \partial}{\xIm \bar\partial}.   
$$ 
On the other hand, we take $w\in L^{n,q}_{(2)}(Y, E\otimes F)\cap \Ker \bar \partial$. 
We put $w^0=\{w_{i_0}\}$, where $w_{i_0}=w|_{U_{i_0}\setminus Z}$. 
We will use $C_i$ to represent some positive constants independent of $w$. 
By Lemma \ref{3.2} below, we have $w^1=\{w^1_{i_0}\}$ such that 
$\bar\partial w^1=w$ on each $U_{i_0}\setminus Z$ with 
$$\|w^1\|^2:=\sum _i\int _{U_i\setminus Z}|w^1_i|^2_{g', h_Eh_F}
\leq C_1\int _{X\setminus Z}|w|^2_{g', h_Eh_F}=C_1\|w\|^2.$$ 
Since $\bar\partial (\delta w^1)=0$, we can obtain $w^2$ such that 
$\bar\partial w^2=\delta w^1$ on each $U_{i_0i_1}\setminus Z$ with 
$$\|w^2\|^2:=\sum _{\{i, j\}\subset I}
\int _{U_{ij}\setminus Z}|w^2_{ij}|^2_{g', 
h_Eh_F}\leq C_2\|w^1\|^2. $$
By repeating this procedure, we obtain $w^q$ such that 
$\bar\partial w^q=\delta w^{q-1}$ with 
$\|w^q\|^2\leq C_q\|w^{q-1}\|^2$. 
In particular, $\|\delta w^q\|^2\leq  C_0\|w\|^2$. 
We put $\beta(w):=\delta w^q=:\{v_{i_0\cdots i_q}\}$. 
Then $\bar\partial v_{i_0\cdots i_q}=0$ and 
$\|v_{i_0\cdots i_q}\|^2<\infty$. 
Thus, $v_{i_0\cdots i_q}\in \Gamma (U_{i_0\cdots i_q}, 
K_X\otimes E\otimes F\otimes \mathcal J(h_F))$ and $\delta (\beta(w))=0$. 
Note that an $E\otimes F$-valued holomorphic 
$(n,0)$-form on $U\setminus Z$, where $U$ is an open subset of $X$,  
with 
a finite $L^2$ norm can be extended 
to an $E\otimes F$-valued holomorphic 
$(n,0)$-form on $U$ (see also Remark \ref{rem215}). 
Therefore, we have a homomorphism 
$$
\bar{\beta}:\frac{L^{n,q}_{(2)}(Y, E\otimes F)\cap \Ker \bar \partial}{\xIm \bar\partial}
\to 
{\check{H}}^q({\mathcal U}, K_X\otimes E\otimes F\otimes \mathcal J(h_F))
$$ 
by the standard diagram chasing. 
It is not difficult to see that $\bar\alpha$ and $\bar\beta$ induce 
the desired isomorphism by the above arguments. 

Next, we note that $\overline {\xIm \bar\partial}=\xIm \bar\partial$ if and 
only if $\overline {\xIm D''^*_{E\otimes 
F}}=\xIm D''^*_{E\otimes F}$ (cf.~\cite[Theorem 1.1.1]{hor}). 
Thus, it is sufficient to prove 
that $\overline {\xIm \bar\partial}=\xIm \bar\partial$. 
Let $w\in \overline {\xIm \bar\partial}$. 
Then there exists a sequence $\{v_k\}\subset 
\xIm \bar\partial$ such that 
$\|w-\bar\partial v_k\|^2\to 0$ if $k\to \infty$. 
By the above construction, $\|\beta (w-\bar\partial v_k)\|^2\to 
0$ when $k\to \infty$. 
This implies that 
$\beta(w-\bar\partial v_k)\to 0$ uniformly on every 
compact subset of $X$. 
Therefore, the image of $w$ 
in $\check {H}^q(\mathcal U, K_X\otimes 
E\otimes F\otimes \mathcal J(h_F))$ is zero because 
$\check {H}^q(\mathcal U, K_X\otimes 
E\otimes F\otimes \mathcal J(h_F))$ is a finite dimensional, 
separated, Fr\'echet space (cf.~\cite[Chap.~VIII, Sec.~A, 19.~Theorem]{gunning-rossi}). 
Thus, $w\in \xIm \bar\partial$ by the above isomorphism. 
For the details of the topology on 
$\mathcal F$ and $H^q(X, \mathcal F)$, where 
$\mathcal F$ is a coherent sheaf on a complex manifold $X$, 
see \cite[\S 55 Coherent Analytic Sheaves as 
Fr\'echet Sheaves]{kk}. 
\end{proof}

There are various formulations for $L^2$-estimates for $\bar\partial$-equations, 
which originated from H\"ormander's paper \cite{hor}. 
The following one is suitable for our purpose. 
We used it in the proof of Claim \ref{cla1}. 

\begin{lem}[$L^2$-estimates for $\bar\partial$-equations on 
complete K\"ahler manifolds]\label{3.2}
Let $U\Subset V$ be small Stein open sets of $X$. If 
$u\in L^{n,q}_{(2)}(U\setminus Z, E\otimes F)_{g', h_Eh_F}$ 
with $\bar\partial u=0$, then 
there exists $v\in L^{n,q-1}_{(2)}(U\setminus Z, E\otimes F)_{g', h_Eh_F}$ such that 
$\bar\partial v =u$. Moreover, 
there exists a positive constant $C$ independent of $u$ 
such that $$\int _{U\setminus Z} |v|^2_{g', h_Eh_F}\leq 
C\int _{U\setminus Z} |u|^2_{g', h_Eh_F}.$$ 
\end{lem}
\begin{proof} 
We can assume that $\omega'=\sqrt{-1}\partial \bar\partial\Psi$ on $V$ because 
$V$ is a small Stein open set. Then 
$(E\otimes F, 
h_E h_Fe^{-\Psi})$ is Nakano positive and $$C_1 h_Eh_F\leq h_Eh_F e^{-\Psi}
\leq C_2 h_Eh_F$$ for some positive constants $C_1$ and 
$C_2$ on $U\setminus Z$. Note that $\Psi$ is a bounded function on $X$ 
by the construction of $g'$. 
It is obvious that $\sqrt{-1}\Theta_{(E\otimes F, h_Eh_Fe^{-\Psi})}\geq_{\Nak} 
\xId_E\otimes \omega'$ on $U\setminus Z$ by the assumption (iii) in 
Theorem \ref{main}. 
Let $w$ be an $E\otimes F$-valued 
$(n,q)$-form on $U\setminus Z$ with measurable 
coefficients. We write 
$$
\| w\|^2=\int _{U\setminus Z} |w|^2_{g', h_Eh_F} dV_{\omega'} \ {\text{ and }}\  
\| w\|^2_0=\int _{U\setminus Z} |w|^2_{g', h_Eh_Fe^{-\Psi}} dV_{\omega'}. 
$$ 
Then $\|w\|$ is finite if and only if $\| w\| _0$ is finite. 
By the well-known $L^2$ estimates for $\bar\partial$-equations 
(cf.~\cite[Th\'eor\`{e}me 4.1, Remarque 4.2]{dem-est} or \cite[(5.1) Theorem]{dem}), 
we obtain an $E\otimes F$-valued $(n, q-1)$-form $v$ on $U\setminus Z$ such 
that $\bar\partial v=u$ and 
$\|v\|^2_0\leq C_0\| u\|^2_0$, where $C_0$ is a positive constant independent of 
$u$. We note that $g'$ is not a complete K\"ahler metric on $U\setminus Z$ but 
$U\setminus Z$ is a complete K\"ahler manifold (cf.~\cite[Th\`eor\'eme 0.2]{dem-est}). 
Therefore, we obtain 
$$
C_1\|v\|^2\leq \|v\|^2_0\leq C_0\|u\|^2_0\leq C_0C_2\|u\|^2. 
$$ 
So, it is sufficient to put $C=\frac{C_0C_2}{C_1}$. 
\end{proof}
Therefore, we obtain 
$$
L^{n,q}_{(2)}(Y,E\otimes F)={\xIm \bar\partial}\oplus 
\mathcal H^{n,q}(E\otimes F)\oplus 
{\xIm D''^*_{E\otimes F}}. 
$$
Thus, $H^q(X, K_X\otimes E\otimes F\otimes \mathcal J(h_F))
\simeq \mathcal H^{n,q}(E\otimes F)$. 

Let $U\Subset V$ be small Stein open sets of $X$. Then there 
exists a smooth strictly psh function $\Phi$ on $V$ such that 
$(L, h_Le^{-\Phi})$ is semi-positive on $V$. 
In this situation, $C'_1\leq e^{-\Phi}\leq C'_2$ on $U$ for some 
positive constants $C'_1$ and $C'_2$. By applying 
the same argument as in Lemma \ref{3.2} to 
$(E\otimes F\otimes L, h_Eh_Fh_Le^{-\Psi-\Phi})$, we obtain 
$$
L^{n,q}_{(2)}(Y,E\otimes F\otimes L)={\xIm \bar\partial}\oplus 
\mathcal H^{n,q}(E\otimes F\otimes L)\oplus 
{\xIm D''^*_{E\otimes F\otimes L}}  
$$ and 
$$H^q(X, K_X\otimes E\otimes F\otimes \mathcal J(h_F)\otimes L)
\simeq \mathcal H^{n,q}(E\otimes F\otimes L)$$ similarly. 

\begin{claim}
The multiplication homomorphism 
$$
\times s:\mathcal H^{n,q}(E\otimes F)\to \mathcal H^{n,q}(E\otimes 
F\otimes L) 
$$ 
is well-defined for every $q\geq 0$. 
\end{claim}
\begin{proof}[Proof of Claim]
By 
Proposition \ref{tsuika}, 
we obtain 
\begin{align*} 
\|D''^*_{E\otimes F}u\|^2+
\|\bar\partial u\|^2=\|{D'}^*u \|^2+
\lla\sqrt{-1}
\Theta(E\otimes F)\Lambda u, u\rra
\end{align*}
for $u\in L^{n,q}_{(2)}(Y, E\otimes F)$, where $\Lambda$ is the adjoint 
of $\omega'\wedge \ \cdot\ $. 
We note that the K\"ahler metric $g'$ on $Y$ is complete. 
If $u\in \mathcal H^{n,q}(E\otimes F)$, 
then 
$D'^*u=0$ and $\langle \sqrt{-1}(\Theta (E)+\xId _E
\otimes \Theta (F))\Lambda u, u\rangle=0$ by the assumption (iii). 
By (iv), we have $\langle \sqrt{-1}(\xId _E\otimes \Theta (L))\Lambda 
u, u\rangle\leq 0$. 
When $u\in \mathcal H^{n,q}(E\otimes F)$, 
$\bar\partial (su)=0$ by the Leibnitz rule and $D'^*(su)=sD'^*u=0$ because 
$s$ is an $L$-valued holomorphic $(0,0)$-form. 
Since $|s|^2_{h_L}$ is a smooth function on $X$, 
there exists a positive number $C$ such that 
$|s|^2_{h_L}<C$ everywhere on $X$. 
Therefore, 
$$
\int _Y |su|^2_{g', h_Lh_Eh_F} dV_{\omega'}<
C\int _Y |u|^2_{g', h_Eh_F}dV_{\omega'}<\infty. 
$$ 
So, $su$ is square integrable. 
We can also see that $su\in \Dom D''^*_{E\otimes F\otimes L}$ 
since $|s|^2_{h_L}<C$ everywhere on $X$. 
Thus, we obtain 
$$
\| D''^*_{E\otimes F\otimes L}(su)\|^2
=\lla\sqrt{-1}\Theta (E\otimes F\otimes L)\Lambda (su), su\rra
$$ 
by Proposition \ref{tsuika}. 
We note that 
\begin{align*}
&\sqrt{-1}(\Theta (E)+\xId _E\otimes \Theta (F)+\xId _E\otimes \Theta (L))
\\
&\geq _{\Nak}(1+\varepsilon)\xId _E\otimes \Theta (L)\\ &\geq _{\Nak}-c'\xId _E\otimes \omega'
\end{align*} 
on $Y=X\setminus Z$ for some constant $c'$. 
On the other hand, 
$$
\langle \sqrt{-1}\Theta (E\otimes F\otimes L)\Lambda (su), su\rangle
=|s|^2\langle \sqrt{-1}
(\xId_E\otimes \Theta (L))\Lambda u, u\rangle \leq 0, 
$$
where $|s|$ is the pointwise norm of $s$ with respect to $h_L$. 
Therefore, $D''^*_{E\otimes F\otimes L}(su)=0$. 
This implies that $su\in \mathcal H^{n,q}(E\otimes F\otimes L)$. 
We finish the proof of the claim. 
\end{proof}
By the above claims, the theorem is obvious because 
$$\times s:\mathcal H ^{n,q}(E\otimes F)\to \mathcal 
H^{n,q}(E\otimes F\otimes L)$$ is injective for every $q$. 
\end{proof}

We close this section with the proof of Corollary \ref{main-cor}. 

\begin{proof}[{Proof of {\em{Corollary \ref{main-cor}}}}]
We put $h_F:=h^{\frac{1}{k}}_D$ as in Example \ref{ex1}, 
where $h_D$ is the natural 
singular hermitian metric on $\mathcal O_X(D)$. 
Then $h_F$ is smooth on $X\setminus D$, $\sqrt {-1}\Theta (F)\geq 0$ in the 
sense of currents, and $\mathcal J(h_F)=\mathcal J(\frac{1}{k}D)$. 
Therefore, we can apply Theorem \ref{main}. 
\end{proof}

\section{Applications:~injectivity and vanishing theorems}\label{sec4}
In this section, we treat only a few applications of Theorem \ref{main}. 
We recommend the reader to see the results in \cite{take} and 
the arguments in \cite[Chapter V, \S 3]{nakayama} for other 
formulations and generalizations. 
See also \cite{fujino2} for various applications 
and generalizations in a more general relative setting. 
For applications in the log minimal model program, which 
can not be covered by the results in this paper, 
see \cite{fujino-non}, \cite{fujino-fundamental}, \cite{fujino-book}, and so on. 

The following formulation is due to Koll\'ar (cf.~\cite[10.13 Theorem]{k2}). 
He stated this result for 
the case where $E$ is a trivial line bundle and $(X,\Delta)$ is klt, 
that is, $\mathcal J(\Delta)\simeq \mathcal O_X$. 

\begin{prop}[Koll\'ar type injectivity theorem]\label{cor03}
Let $f:X\to Y$ be a proper surjective morphism from a 
compact K\"ahler manifold $X$ to a normal projective variety $Y$. 
Let $L$ be a holomorphic line bundle on $X$ and 
let $D$ be an effective divisor on $X$ such that $f(D)\ne Y$. 
Assume that $L\equiv f^*M+\Delta$, where $M$ is a nef and 
big $\mathbb Q$-divisor on $Y$ and $\Delta$ is an effective $\mathbb Q$-divisor 
on $X$. 
Let $(E,h_E)$ be a Nakano 
semi-positive holomorphic vector bundle on $X$. 
Then 
$$
H^q(X, K_X\otimes E\otimes L\otimes \mathcal J(\Delta))\to 
H^q(X, K_X\otimes E\otimes L\otimes \mathcal O_X(D)\otimes \mathcal J(\Delta))
$$ 
is injective for every $q\geq 0$, 
where $\mathcal J(\Delta)$ is the multiplier ideal sheaf 
associated to the effective $\mathbb Q$-divisor $\Delta$. 
\end{prop}
\begin{proof}
By taking $P\in \Pic^0(X)$ suitably, we have $L\otimes P\sim_{\mathbb Q}
f^*M+\Delta$. 
We can assume that 
$L\sim_{\mathbb Q}f^*M+\Delta$ by replacing $L$ (resp.~$E$) 
with $L\otimes P$ (resp.~$E\otimes P^{-1}$). 
By Kodaira's lemma (see \cite[Proposition 2.61]{km}), 
we can further assume that $M$ is ample (cf.~Definition \ref{28}). 
Let $h:=\Phi_{|mM|}:Y\to \mathbb P^N$ be the 
embedding induced by the complete linear system $|mM|$ 
for a large integer $m$. 
Then $\mathcal O_Y(mM)\simeq h^*\mathcal O_{\mathbb P^N}(1)$. 
We can take an effective divisor $A$ on $\mathbb P^N$ such 
that $\mathcal O_{\mathbb P^N}(A)\simeq 
\mathcal O_{\mathbb P^N}(l)$ for some 
positive integer $l$ 
and $D'=f^*h^*A-D$ is an effective divisor on $X$. 
We add $D'$ to $D$ and 
can assume that $D=f^*h^*A$. 
Under these extra assumptions, we can easily construct hermitian 
metrics satisfying the assumptions in 
Theorem \ref{main} (see Example \ref{ex1}). We finish the proof of the 
proposition. 
\end{proof}
\begin{rem}[Numerical equivalence] 
In the above proposition, we note that $L\equiv f^*M+\Delta$ 
means $c_1(L)=c_1(f^*M+\Delta)$ in $H^2(X, \mathbb R)$, 
where $c_1$ is the first Chern class of $\mathbb Q$-divisors or 
line bundles. 
\end{rem}
\begin{rem}
Proposition \ref{cor03} is a generalization of \cite[10.13 Theorem]{k2}, 
which is stated for a compact K\"ahler manifold. 
However, the proof of \cite[10.13 Theorem]{k2} given in \cite{k2} 
works only for {\em{projective manifolds}}. 
In \cite[10.17.3 Claim]{k2}, we need an ample divisor on $X$ to 
prove local vanishing theorems. 
\end{rem}

The following proposition is a reformulation 
of \cite[5.12.~Corollary 
b)]{ev} from the analytic viewpoint. 
It is essentially the same as Proposition \ref{cor03}. 
In \cite{ev}, $E$ is trivial and $\mathcal J\simeq \mathcal O_X$. 

\begin{prop}[Esnault--Viehweg type injectivity theorem]\label{ev-type} 
Let $X$ be a smooth projective variety and 
let $D$ be an effective divisor on $X$. 
Let $(E, h_E)$ be a Nakano semi-positive holomorphic 
vector bundle and let $L$ be a holomorphic 
line bundle on $X$. 
Assume that $L^{\otimes k}(-D)$ is nef and abundant, that is, 
$\kappa (L^{\otimes k}(-D))=\nu(L^{\otimes k}(-D))$, for 
some positive integer $k$. 
Let $B$ be an effective divisor on $X$ such that 
$$H^0(X, (L^{\otimes k}(-D))^{\otimes l}\otimes 
\mathcal O_X(-B))\ne 0$$ for 
some $l>0$. 
Then 
$$
H^q(X, K_X\otimes E\otimes L\otimes \mathcal J)
\to H^q(X, K_X\otimes E\otimes L\otimes 
\mathcal J\otimes \mathcal O_X(B))
$$ 
is injective for every $q$, where 
$\mathcal J:=\mathcal J(\frac{1}{k}D)$ is 
the multiplier ideal sheaf associated to 
the effective $\mathbb Q$-divisor 
$\frac{1}{k}D$. 
\end{prop}
\begin{proof}
Let $\pi:Z\to X$ be a projective birational morphism from 
a smooth projective variety $Z$ with the following 
properties:~(i) 
There exists a proper surjective morphism between 
smooth projective varieties $f:Z\to Y$ with 
connected fibers, and (ii) there is a nef and big 
$\mathbb Q$-divisor $M$ on $Y$ such 
that $\pi^*(L^{\otimes k}(-D))\sim_{\mathbb Q}f^*M$. 
For the proof, see \cite[Proposition 2.1]{ka}. 
On the other hand, $R^i\pi_*(K_{Z/X}\otimes \mathcal J(\frac{1}{k}\pi^*D))=0$ 
for $i>0$ and 
$\pi_*(K_Z\otimes \mathcal J(\frac{1}{k}\pi^*D))\simeq 
K_X\otimes\mathcal J(\frac{1}{k}D)$ by 
\cite[Theorem 9.2.33, and Example 9.6.4]{laz}. 
We note that $(\pi^*E, \pi^*h_E)$ is Nakano semi-positive on $Z$. 
So, we can assume that $X=Z$ without loss of generality. 
It is not difficult to see that $f(B)\neq Y$ by the assumption 
that 
$H^0(X, (L^{\otimes k}(-D))^{\otimes l}\otimes \mathcal O_X(-B))\ne 0$ for 
some $l>0$. 
Thus, this proposition follows from Proposition \ref{cor03}. 
\end{proof}

The referee pointed out that 
Proposition \ref{ev-type} is sharper than \cite[Theorem 3.1]{ep}. 

\begin{rem}In this remark, we use the notation in \cite[Theorem 3.1]{ep}. 
Let $k$ be a positive integer 
such that $k>\lambda$. 
We take general members $D_1, \cdots, D_k$ of $H^0(X, A\otimes \mathfrak a)$ and 
put $$
D=\frac{\lambda}{k}(D_1+\cdots +D_k). 
$$ 
Then $L-D$ is nef and abundant and $L-D-\epsilon B$ 
is $\mathbb Q$-effective for some $0<\epsilon <1$. 
By the construction, we have 
$\mathcal J(D)=\mathcal J(X, \mathfrak a^{\lambda})$ (cf.~\cite[Proposition 9.2.28]
{laz}). Therefore, by Proposition \ref{ev-type} for 
$E=\mathcal O_X$, we obtain that 
$$
H^i(X, \mathcal O_X(K_X+L)\otimes \mathcal J(X, \mathfrak a^{\lambda}))
\to H^i(X, \mathcal O_X(K_X+L+B)\otimes 
\mathcal J(X, \mathfrak a^{\lambda}))
$$ 
is injective for every $i$.  
\end{rem}

\begin{rem}[Vanishing theorem and torsion-freeness] 
Proposition \ref{cor03} gives some 
generalizations of Koll\'ar's vanishing and torsion-free theorems. 
We do not pursue them here because 
we discuss them in a more general relative setting 
in \cite{fujino2}. We just mention that \cite[10.15 Corollary]{k2} 
holds for $K_X\otimes E\otimes \mathcal J(\Delta)$, where we use 
the same notation as in Proposition \ref{cor03}. 
We note \cite[Example 9.5.9]{laz} when we restrict the multiplier ideal sheaf 
$\mathcal J(\Delta)$ to a {\em{general}} hypersurface. 
Related topics are in \cite[6.12 Corollary, and 6.17 Corollary]{ev} and 
\cite[Section 3]{ep}. 
\end{rem}

By combining Proposition \ref{cor03} with 
Serre's vanishing theorem, we obtain the next corollary. 
It may be better to 
be called Nadel type vanishing theorem. 

\begin{cor}[Kawamata--Viehweg type vanishing theorem]\label{cor2}  
Let $X$ be a smooth projective variety and let $L$ be a holomorphic line bundle 
on $X$. 
Assume that $L\equiv M+\Delta$, where $M$ is a nef and big 
$\mathbb Q$-divisor on $X$ and 
$\Delta$ is an effective $\mathbb Q$-divisor on $X$. 
Let $(E, h_E)$ be a Nakano semi-positive 
holomorphic vector bundle on $X$. 
Then $H^q(X, K_X\otimes E\otimes L\otimes \mathcal J(\Delta))=0$ 
for $q\geq 1$. 
Moreover, 
if $\Delta$ is integrable outside finitely many points, 
then $H^q(X, K_X\otimes E\otimes L)=0$ 
for $q\geq 1$. 
\end{cor}
\begin{proof}
We use Proposition \ref{cor03} 
under the assumption that $Y = X$ and
$f=\xid _X$. We take an effective ample divisor $D$ on $X$ and apply Proposition 
\ref{cor03}. Then we obtain that
$$
H^q(X, K_X\otimes E\otimes L\otimes \mathcal J(\Delta))\to 
H^q(X, K_X\otimes E\otimes L\otimes \mathcal J(\Delta)\otimes 
\mathcal O_X(mD))
$$
is injective for $m>0$ and $q \geq  0$. 
By Serre's vanishing theorem, we
have $H^q(X,K_X \otimes E \otimes  L \otimes  
\mathcal J(\Delta)) = 0$ for $q\geq 1$. 
When $\Delta$ is integrable outside finitely many points, $\mathcal O_X/\mathcal J(\Delta)$ 
is a skyscraper sheaf. Therefore,
$H^q (X, K_X\otimes  E \otimes L \otimes \mathcal O_X/\mathcal J(\Delta))= 0$ 
for $q\geq 1$. 
By combining it with the above mentioned vanishing result, we obtain 
the desired result. 
\end{proof}

The final result is a slight generalization of Demailly's 
formulation of Kawamata--Viehweg type vanishing theorem. 

\begin{cor}[{cf.~\cite[Main Theorem]{trans}}]\label{sai}
Let $L$ be a holomorphic line bundle on an $n$-dimensional 
projective manifold $X$. 
Assume that some positive power $L^{\otimes k}$ can be written $L^{\otimes 
k}\simeq M\otimes \mathcal O_X(D)$, where 
$M$ is a nef line bundle and $D$ is an 
effective divisor such that 
$\frac{1}{k}D$ is integrable on $X\setminus B$. 
Let $\nu=\nu(M)$ be the numerical dimension of the nef line bundle $M$. 
Let $(E, h_E)$ be a Nakano semi-positive holomorphic vector bundle on $X$. 
Then 
$H^q(X, K_X\otimes E\otimes L)=0$ for 
$q>n-\min\{\max\{\nu, \kappa (L)\}, \codim B\}$. 
\end{cor}

\begin{proof}[Sketch of the proof] 
By the standard slicing arguments, 
we can reduce it to the case
where $\min\{\max\{\nu, \kappa (L)\}, \codim B\}=\dim X$. 
We note that $\codim B=\infty$ if and only if $B=\emptyset$. 
By Kodaira's 
lemma (cf.~\cite[Lemma 2.60, Proposition 2.61]{km}), 
we can further reduce it to the case when 
$M$ is ample. 
We note that if $A$ is a general smooth very ample 
Cartier divisor on $X$ then 
\begin{align*}
0\to K_X\otimes E\otimes L\to K_X\otimes E\otimes L\otimes \mathcal O_X(A)
\to K_A\otimes E|_A\otimes L|_A\to 0
\end{align*} 
is exact and 
$\mathcal J(\frac{1}{k}D)|_A=\mathcal J(\frac{1}{k}D|_A)$. 
In particular, $\frac{1}{k}D|_A$ is integrable on $A\setminus B|_A$. 
For the details of these reduction
arguments, see the first and second steps in the proof of the main 
theorem in \cite{trans}. 
Therefore, this corollary follows from the previous
corollary:~Corollary \ref{cor2}.
\end{proof}

\ifx\undefined\bysame
\newcommand{\bysame|{leavemode\hbox to3em{\hrulefill}\,}
\fi


\begin{thebibliography}{FST}

\bibitem[D1]{dem-est} 
J-P.~Demailly, 
Estimations $L\sp{2}$ pour l'op\'erateur $\bar \partial $ d'un fibr\'e 
vectoriel holomorphe semi-positif au-dessus d'une vari\'et\'e k\"ahl\'erienne 
compl\`{e}te, Ann. Sci. \'Ecole Norm. Sup. (4)  {\textbf{15}} 
(1982),  no. 3, 457--511. 

\bibitem[D2]{dem-coh} 
J-P.~Demailly, 
Cohomology of $q$-convex spaces in top 
degrees, 
Math. Z. {\textbf{204}} 
(1990), no. 2, 283--295.

\bibitem[D3]{trans} 
J-P.~Demailly, 
Transcendental proof of a generalized 
Kawamata--Viehweg vanishing  theorem, 
Geometrical and algebraical aspects in 
several complex variables  (Cetraro, 1989),  81--94, 
Sem. Conf., {\textbf{8}}, EditEl, Rende,  1991. 
 
\bibitem[D4]{dem} 
J-P.~Demailly, 
Multiplier ideal sheaves and analytic 
methods in algebraic geometry, 
School on Vanishing Theorems and Effective 
Results in Algebraic Geometry 
(Trieste, 2000), 
1--148, ICTP Lect. Notes, {\textbf{6}}, Abdus 
Salam Int. Cent. Theoret. Phys., Trieste,  2001.

\bibitem[EP]{ep} 
L.~Ein, M.~Popa, 
Global division of cohomology classes via 
injectivity, 
Special volume in honor of 
Melvin Hochster. Michigan Math. J. {\textbf{57}} (2008), 249--259.

\bibitem[E]{enoki}
I.~Enoki, 
Kawamata--Viehweg vanishing theorem for compact K\"ahler 
manifolds, 
Einstein metrics and Yang--Mills connections 
(Sanda, 1990), 
59--68, Lecture Notes in Pure and Appl. 
Math., {\textbf{145}}, 
Dekker, New York, 1993. 

\bibitem[EV]{ev}
H.~Esnault, E.~Viehweg, 
{\em{Lectures on vanishing theorems}}, 
DMV Seminar, {\textbf{20}}. Birkh\"auser Verlag, Basel, 1992. 

\bibitem[F1]{fujino2}O.~Fujino, 
A transcendental approach to Koll\'ar's injectivity theorem II, 
preprint (2006). 

\bibitem[F2]{fujino3}O.~Fujino, 
On Koll\'ar's injectivity theorem (Japanese), 
S\=urikaisekikenky\=usho K\=oky\=uroku 
No. 1550 (2007), 131--140. 

\bibitem[F3]{fujino-multi}O.~Fujino, 
Multiplication maps and vanishing theorems for toric varieties, 
Math. Z. {\textbf{257}} (2007), no. 3, 631--641.

\bibitem[F4]{fujino-vani} 
O.~Fujino, 
Vanishing theorems for toric polyhedra, 
Higher dimensional algebraic varieties and vector bundles, 
81--95, RIMS K\^oky\^uroku Bessatsu, B9, Res. Inst. Math. Sci. 
(RIMS), Kyoto, 2008.

\bibitem[F5]{fujino-on} 
O.~Fujino, 
On injectivity, vanishing and torsion-free theorems for 
algebraic varieties, 
Proc. Japan Acad. Ser. A Math. Sci. {\textbf{85}} (2009), no. 8, 95--100.

\bibitem[F6]{fujino-eff1} 
O.~Fujino, 
Effective base point free theorem for log canonical 
pairs---Koll\'ar type theorem, 
Tohoku Math. J. (2) {\textbf{61}} (2009), no. 4, 475--481. 

\bibitem[F7]{fujino-eff2} 
O.~Fujino, 
Effective base point free theorem for log canonical pairs, 
II.~Angehrn--Siu type theorems, 
Michigan Math. J. {\textbf{59}} (2010), no. 2, 303--312.

\bibitem[F8]{fujino-theory} 
O.~Fujino, 
Theory of non-lc ideal sheaves: basic properties, 
Kyoto J. Math. {\textbf{50}} (2010), no. 2, 225--245. 

\bibitem[F9]{fujino-non} 
O.~Fujino, 
Non-vanishing theorem for log canonical pairs, 
to appear in Journal of Algebraic Geometry.

\bibitem[F10]{fujino-qlog} 
O.~Fujino, 
Introduction to the theory of quasi-log varieties, 
{\em{Classification of Algebraic Varieties}}, 
289--303, EMS Ser. of Congr. Rep. Eur. Math. Soc., Z\"urich, 2010.

\bibitem[F11]{fujino-fundamental} 
O.~Fujino, 
Fundamental theorems for the log minimal model program, 
to appear in Publ. Res. Inst. Math. Sci.

\bibitem[F12]{fujino-book} 
O.~Fujino, 
Introduction to the log minimal model program for log canonical pairs, 
preprint (2009). arXiv:0907.1506v1

\bibitem[F13]{fujino-surface} 
O.~Fujino, 
Minimal model theory for log surfaces, preprint (2010).

\bibitem[F14]{fujino-bpf} 
O.~Fujino, Base point free theorems---saturation, b-divisors, and 
canonical bundle formula---, preprint (2011). 

\bibitem[FST]{fst} 
O.~Fujino, K.~Schwede, and S.~Takagi, 
Supplements to non-lc ideal sheaves, 
to appear in RIMS K\^oky\^uroku Bessatsu. 

\bibitem[Fk]{fukuda} 
S.~Fukuda, 
On numerically effective log canonical divisors, 
Int. J. Math. Math. Sci. {\textbf{30}} (2002), no. 9, 521--531.

\bibitem[GR]{gunning-rossi} 
R.~C.~Gunning, H.~Rossi, 
{\em{Analytic functions of several complex variables}}, 
Prentice-Hall, Inc., Englewood Cliffs, N.J. 1965. 

\bibitem[H]{hor}
L.~H\"ormander, 
$L\sp{2}$ estimates 
and existence theorems for the 
$\bar \partial $  operator, Acta Math. {\textbf{113}} (1965), 89--152.

\bibitem[KK]{kk} 
L.~Kaup, B.~Kaup, 
{\em{Holomorphic functions of several variables}}, 
{\em{An introduction to the fundamental theory}}, 
With the assistance of Gottfried Barthel. Translated 
from the German by Michael Bridgland. 
de Gruyter Studies in Mathematics, {\textbf{3}}. Walter de Gruyter $\&$ Co., Berlin, 1983.

\bibitem[Ka]{ka} 
Y.~Kawamata, 
Pluricanonical systems on minimal algebraic varieties, 
Invent. Math. {\textbf{79}} (1985),  no. 3, 567--588. 

\bibitem[Ko1]{k1} 
J.~Koll\'ar, 
Higher direct images of dualizing sheaves. I, 
Ann. of Math. (2) {\textbf{123}} (1986), no. 1, 11--42.

\bibitem[Ko2]{k2} 
J.~Koll\'ar, {\em{Shafarevich maps and automorphic forms}}, 
M. B. Porter Lectures. Princeton University Press, 
Princeton, NJ, 1995. 

\bibitem[KM]{km} 
J.~Koll\'ar, S.~Mori, {\em{Birational geometry of algebraic varieties}}. 
With the collaboration of C. H. Clemens and A. Corti. 
Translated from the 1998 Japanese original. 
Cambridge Tracts in Mathematics, {\textbf{134}}. Cambridge University 
Press, Cambridge, 
1998. 

\bibitem[L]{laz} 
R.~Lazarsfeld, {\em{Positivity in algebraic geometry. II. Positivity for vector 
bundles, and multiplier ideals}}, 
Ergebnisse der Mathematik und ihrer Grenzgebiete. 3. Folge. A 
Series of  Modern Surveys in Mathematics 
[Results in Mathematics and Related Areas. 3rd  Series. A 
Series of Modern Surveys in Mathematics], {\textbf{49}}. 
Springer-Verlag, Berlin,  2004. 

\bibitem[Nd]{nadel} 
A.~Nadel, 
Multiplier ideal sheaves and existence of K\"ahler-Einstein metrics of 
positive scalar curvature, 
Proc. Nat. Acad. Sci. U.S.A. {\textbf{86}} (1989), no. 19, 7299--7300.

\bibitem[Nk]{nakayama}
N.~Nakayama, 
{\em{Zariski-decomposition and abundance}}, 
MSJ Memoirs, {\textbf{14}}. 
Mathematical Society of Japan, Tokyo, 2004. 

\bibitem[O1]{oh-hon}
T.~Ohsawa, 
{\em{Analysis of several complex variables}}, 
Translated from the Japanese by Shu Gilbert Nakamura. 
Translations of Mathematical Monographs, {\textbf{211}}. 
Iwanami Series in Modern Mathematics. 
American Mathematical Society, Providence, RI,  
2002. 

\bibitem[O2]{ohsawa} 
T.~Ohsawa, 
On a curvature condition that implies a cohomology injectivity 
theorem of Koll\'ar--Skoda type, 
Publ. Res. Inst. Math. Sci. {\textbf{41}} (2005), no. 3, 565--577.

\bibitem[S]{siu} 
Y.-T.~Siu, 
A vanishing theorem for semipositive line bundles over non-K\"ahler manifolds, 
J. Differential Geom. {\textbf{19}} (1984), no. 2, 431--452. 
 
\bibitem[Tk]{take}
K.~Takegoshi, 
Higher direct images of canonical sheaves tensorized with 
semi-positive vector bundles by proper K\"ahler morphisms, 
Math. Ann.  {\textbf{303}} (1995), no. 3, 389--416. 

\bibitem[Tn]{tan} 
S.~G.~Tankeev, On 
$n$-dimensional canonically polarized varieties 
and varieties of fundamental type, 
Math. USSR-Izv. vol. {\textbf{5}} (1971), no. 1, 29--43.  
\end{thebibliography}
\end{document}